\documentclass[11pt,reqno,oneside]{amsart}
\usepackage[asymmetric,top=3.5cm,bottom=4.0cm,left=3.1cm,right=3.1cm]{geometry}
\usepackage{amsmath,amsfonts,amsthm,mathrsfs,amssymb, cite,array}
\usepackage[usenames]{color}
\usepackage{mathtools} 
\usepackage{framed}
\usepackage{enumerate}
\usepackage{graphicx}      
\usepackage{caption}
\usepackage{subcaption}

\usepackage{hyperref}
\usepackage{xcolor}
\hypersetup{
  colorlinks,
  linkcolor={blue!80!black},
  urlcolor={blue!80!black},
  citecolor={blue!80!black}
}
\usepackage[capitalize,noabbrev]{cleveref}
\usepackage{graphicx}
\usepackage{latexsym}
\usepackage{enumerate}
\usepackage{booktabs}
\usepackage{adjustbox}
\usepackage{multirow}

\newtheorem{thm}{Theorem}[section]

\newtheorem{Definition}[thm]{Definition}
\numberwithin{equation}{section}
\newtheorem{remark}[thm]{Remark}

\newcommand{\Bu}{\boldsymbol{u}}

\author{}
\date{\today}
\title[Wave Front Tracking for 2×2 systems]{A Wave Front Tracking Scheme for Flux Reconstruction in $2\times 2$ Hyperbolic Conservation Laws}

\author{Chaohua Duan}
\address{Department of Mathematics, City University of Hong Kong, Hong Kong SAR, China}
\email{chduan3-c@my.cityu.edu.hk, chduan3@gmail.com}

\author{Yan Jiang}
\address{Department of Mathematics, City University of Hong Kong, Hong Kong SAR, China}
\email{yjian24@cityu.edu.hk}

\author{Hongyu Liu}
\address{Department of Mathematics, City University of Hong Kong, Hong Kong SAR, China}
\email{hongyu.liuip@gmail.com, hongyliu@cityu.edu.hk}

\author{Wenjian Peng}
\address{Department of Mathematics, City University of Hong Kong, Hong Kong SAR, China}
\email{wjpeng5-c@my.cityu.edu.hk}

\allowdisplaybreaks

\begin{document}

\begin{abstract}
This paper introduces a novel wave front tracking framework for reconstructing unknown flux functions in $2\times 2$ hyperbolic conservation laws, extending beyond the well-studied scalar case. By analyzing Riemann solutions at fixed observation times, we develop explicit reconstruction formulas that handle arbitrary combinations of shock and rarefaction waves through a unified equivalent shock concept. Our method constructs piecewise quadratic $C^1$ flux approximations with rigorous convergence guarantees: the approximation errors decrease quadratically with the discretization parameters for function values and linearly for derivatives under $C^{1,1}$ regularity, with enhanced cubic and quadratic convergence respectively under $C^3$ regularity. Applications to the isentropic Euler equations and the mathematically equivalent p-system in compressible fluid dynamics demonstrate the method's capability to identify complete equations of state from limited dynamic measurements, providing a systematic approach to a fundamental inverse problem in continuum mechanics.

\medskip

    
\end{abstract}

   \maketitle   
    
\section{Introduction}

\subsection{Background}

Hyperbolic systems of conservation laws arise naturally in numerous physical contexts, including fluid dynamics, gas dynamics, electromagnetism, and traffic flow. Consider the following hyperbolic system of conservation laws:
\begin{equation}\label{main:eq0}
\begin{cases}
\partial_t \Bu(t, x) + \operatorname{div}_x\left(F(\Bu(t, x))\right) =0,\\
\Bu(0,x) = \bar{\Bu}(x), 
\end{cases}
\end{equation}
where $\Bu: (t, x)\in\mathbb{R}_+\times\mathbb{R}^{n}\rightarrow \mathbb{R}^m$ is the vector of conserved quantities (state variables), and $F: \mathbb{R}^m \rightarrow \mathbb{R}^{m\times n}$ is the \emph{flux function}. Here, $m\geq 1$ denotes the number of conservation laws, $t\in\mathbb{R}_+$ and $x\in\mathbb{R}^n$  respectively signify the time and space variables. A central object governing the dynamics of such systems is the flux function $F(\Bu)$, which encodes the constitutive relation between the conserved variables and their transport mechanisms. Determining the precise form of this flux function is of fundamental importance, as it directly dictates the propagation of nonlinear waves-such as shocks, rarefactions, and contact discontinuities-and hence controls both the qualitative and quantitative features of the system's evolution.

However, in many real-world applications, the explicit expression of $F$ is unknown or only partially accessible. This situation often arises in complex or composite materials, where the constitutive law between stress and strain is not explicitly measurable; in plasma physics and astrophysical flows, where the flux depends on intricate electromagnetic or thermodynamic interactions; in biological and physiological systems, where macroscopic transport laws emerge from poorly understood microscopic processes; and even in traffic or crowd dynamics, where the flux-density relationship (the so-called fundamental diagram) must be inferred from empirical data. In each of these settings, direct measurement of $F$ is either impractical or impossible, making its identification a problem of both theoretical and practical significance.

Consequently, the \emph{inverse problem} of determining the underlying flux function from observational data represents more than a technical challenge-it provides a fundamental pathway to discovering the intrinsic physical laws governing complex systems. Successfully reconstructing $F$ from solution data implies that we can deduce the constitutive relations solely by observing the system's dynamic response, effectively bridging the gap between emergent macroscopic phenomena and their microscopic underpinnings. This approach transforms the flux function from a postulated model into an empirically determined quantity, offering a mathematically rigorous framework for model discovery and validation.

In this paper, we address this foundational challenge by developing a systematic methodology for flux inversion. Our approach relies on multiple active measurements and the analysis of Riemann problems, which form the fundamental building blocks for understanding general solutions to hyperbolic conservation laws. As observed by Courant and Friedrichs \cite{CM2014}, Riemann solutions exhibit characteristic wave patterns-either continuous rarefaction waves or discontinuous shock waves-that encode essential information about the underlying flux structure.

\subsection{Technical developments and discussion}

Hyperbolic systems of conservation laws constitute one of the most fundamental classes of partial differential equations in mathematical physics, providing the theoretical foundation for describing wave propagation phenomena across diverse physical domains. These systems arise naturally in gas dynamics, astrophysics, elasticity theory, electromagnetism, and traffic flow modeling, capturing essential features of wave interactions, shock formation, and energy propagation in nonlinear media. The intrinsic finite-speed propagation characteristics and the emergence of discontinuous solutions make hyperbolic conservation laws particularly suited for modeling physical systems where wave phenomena dominate.

The mathematical theory of hyperbolic conservation laws has developed significantly since Riemann's pioneering work in 1860 \cite{Riemann1860}. In 1883, Stokes \cite{Stokes1883} studied the formation of singularities for a class of scalar conservation laws. For systems of conservation laws, Lax \cite{Lax1957} provided in 1957 a rigorous definition of different types of discontinuities and established the existence of solutions to Riemann problems for small discontinuities. In 1964, Lax \cite{Lax1964} introduced a diagonalization procedure for $2\times 2$ systems. The following year, Glimm \cite{Glimm1965} used the Glimm scheme to prove the existence of BV solutions for initial data with small bounded variation. Subsequently, DiPerna \cite{Diperna1967} employed the wave-front tracking method and proved the existence of solutions for genuinely nonlinear systems. The front-tracking methodology was further developed by Risebro \cite{Risebro1993}, who provided an alternative to the random choice method, and systematically presented in the monograph by Holden and Risebro \cite{HoldenRisebro2015}. In 1992, Bressan \cite{Bressan1992} extended these results to $n\times n$ systems that are not necessarily genuinely nonlinear. These foundational works have established a robust mathematical framework for understanding solution behavior in specific settings, including the well-posedness theory developed by Bressan and colleagues \cite{BA1995, BA2000, BSCR} for one-dimensional systems with small variation. Readers may also refer to \cite{Lax1960, HN2003, Da2005, MA1984, Serre1999, Serre2000} for comprehensive treatments.

Substantial challenges emerge beyond the $2\times 2$ case. When the number of unknown functions exceeds two, the diagonalization method may fail. The interactions among different characteristic families become more complex and harder to estimate. Furthermore, in higher spatial dimensions, systems of conservation laws introduce new phenomena and profound challenges. Since characteristic lines generalize to characteristic surfaces in higher dimensions, one must account for numerous complex situations including the coupling of multiple characteristic fields, intricate wave interactions (such as shock-shock collisions), and the complex geometry of wave curves in phase space. These difficulties are compounded when considering systems in higher spatial dimensions, where even the forward problem theory remains markedly less developed, facing obstacles such as the complex geometry of wave fronts, the instability of planar shocks, and the absence of a comprehensive well-posedness theory for general data.

The development of the forward theory exposed the profound mathematical challenges inherent in these systems. The inverse problem, however, remains largely uncharted territory, as it must not only overcome these established hurdles but also grapple with the additional layer of uncertainty introduced by having to identify, rather than simply evolve, the underlying flux laws. Existing results are predominantly confined to one-dimensional scalar equations, with the literature for systems being exceptionally sparse-a direct reflection of the fact that the forward theory for systems is itself significantly more complex. The primary challenge stems from and is compounded by the fundamental nonlinearity of characteristic structures and the sensitive dependence of solutions on flux properties, which are defining features of the forward problem. Previous approaches have been largely limited to simplified settings that bypass these core difficulties, including optimization-based strategies for specific model forms \cite{BCS2009}, stability analysis via Carleman estimates \cite{FT2021}, reconstruction methods building upon front-tracking algorithms \cite{Risebro1993, HoldenRisebro2015, HPR2014}, and flux identification from shock asymptotics in one-dimensional scalar cases \cite{KT2005}. A notable progress was recently made in \cite{DLMW2025}, which introduced an operator-theoretic framework for determining nonlinear balance laws in product space using boundary measurements-marking a systematic approach to a broader class of inverse problems.

 The flux inversion problem for general coupled $2\times 2$ hyperbolic systems presents unique and formidable challenges that not only mirror but magnify the difficulties encountered in the forward problem. Just as the analysis of the forward problem for systems must contend with the intricate coupling between characteristic fields and complex interactions of multiple wave families, the inverse problem must now use limited observational data to uniquely disentangle the simultaneous contributions of multiple unknown flux functions from these very same complex wave patterns. This multi-parameter inversion is fundamentally more ill-posed than its scalar counterpart. Furthermore, the nonlinear dependence of wave speeds on flux derivatives-a key feature that governs shock formation and wave interactions in the forward problem-becomes a severe complicating factor in the inverse setting, rendering traditional methods that rely on dense temporal data or strong prior assumptions on functional forms inadequate. Addressing this challenge thus requires developing new methodologies that directly engage with the core mathematical structure of hyperbolic systems, representing a critical step forward in the theory of inverse problems for conservation laws.

In this work, we address the flux inversion problem for general coupled hyperbolic systems of the form \eqref{main:eq0} through a novel algorithmic reconstruction framework. Our method employs an iterative procedure using multiple carefully designed Riemann-type initial conditions to probe the system's response. Specifically, we construct a sequence of Riemann problems with progressively refined initial states $\Bu_h = (u_h, v_h)$ for $h = 0, \ldots, 2^m$, where the grid spacings $\delta$ and $\eta$ decrease exponentially with refinement level $m$. For each Riemann problem, we observe the solution at a fixed time $T$ and extract flux information by analyzing the resulting wave patterns-whether shock waves, rarefaction waves, or their combinations. This multi-observation strategy enables comprehensive flux reconstruction over the entire parameter range of interest through a systematic accumulation of local flux information.

The significance of our results extends to fundamental applications in continuum mechanics, particularly in compressible fluid dynamics where our method enables the determination of equations of state $p(\rho)$ through multiple carefully designed dynamic experiments. The unified treatment of the isentropic Euler equations and its mathematically equivalent p-system formulation demonstrates the method's versatility across different physical descriptions. This systematic multi-experiment approach proves especially valuable for characterizing complex material behavior, as it reconstructs complete constitutive relations from wave propagation patterns observed in varied initial conditions.

Looking forward, the algorithmic framework developed in this work opens new avenues for experimental determination of constitutive relations in continuum mechanics. By bridging the gap between mathematical theory and practical implementation, our approach provides both a rigorous foundation and a practical methodology for flux identification in hyperbolic systems-a fundamental challenge that has remained largely open despite decades of research in conservation laws theory. The iterative nature of our reconstruction algorithm, combined with its provable convergence properties, makes it suitable for implementation in computational inverse problems where flux functions need to be determined from experimental data.

The remainder of this paper is structured as follows. Section \ref{Sec:main-result} presents the mathematical framework and main theoretical results. Section \ref{Sec:proof} provides detailed proofs of the flux reconstruction methodology across all wave interaction patterns. Section \ref{Sec:Application} demonstrates applications to physically significant systems, validating the method's practical utility across different domains of continuum mechanics.

\section{Preliminaries and Main Results}\label{Sec:main-result}

\subsection{Mathematical Framework}

We primarily focus our investigation on the following general coupled system of conservation laws in one space dimension:
\begin{equation}\label{Main:eqf}
  \begin{cases}
    \partial_t u + \partial_x f_{1}(v)=0, \\
    \partial_t v+\partial_x f_{2}(u)=0, \\
    u(0,x) = u_0(x), \quad v(0,x) = v_0(x).
  \end{cases}
\end{equation}
Here, $u(t, x), v(t, x): \mathbb{R}_+\times\mathbb{R}\rightarrow \mathbb{R}$ are the conserved quantities (state variables), and $f_1, f_2: \mathbb{R} \rightarrow \mathbb{R}$ represent the unknown flux functions to be reconstructed. We assume sufficient smoothness to ensure the existence of a Standard Riemann Semigroup (SRS), see Definition \ref{def:SRS}. This specific coupling structure, where each variable's evolution depends on the flux of the other variable, captures essential features of many physical systems while maintaining mathematical tractability for our inversion methodology.

The system can be compactly expressed in vector form as:
\begin{equation}
\partial_t \Bu + \partial_x F(\Bu) = 0, \quad \text{with } \Bu = \begin{pmatrix} u \\ v \end{pmatrix}, \quad F(\Bu) = \begin{pmatrix} f_1(v) \\ f_2(u) \end{pmatrix}.
\end{equation}

The well-posedness theory for hyperbolic conservation laws is rigorously established through the framework of Standard Riemann Semigroups \cite{BA1995, BA2000}. We begin by recalling the fundamental definition:

\begin{Definition}[Standard Riemann Semigroup]\label{def:SRS}
Let $\Omega \subset \mathbb{R}^2$ be an open domain. A flux function $F: \Omega \to \mathbb{R}^2$ belongs to $\operatorname{Hyp}(\Omega)$ if $F$ is sufficiently smooth and strictly hyperbolic in $\Omega$, i.e., the Jacobian $DF(\Bu)$ has real distinct eigenvalues for all $\Bu \in \Omega$. 

A mapping $S^F: [0,+\infty) \times \mathscr{D}^F \to \mathscr{D}^F$ is called a Standard Riemann Semigroup for the system \eqref{Main:eqf} if $F \in \operatorname{Hyp}(\Omega)$ and the following properties hold:
\begin{itemize}
    \item \textbf{Semigroup structure}: $S_0^F \Bu = \Bu$ and $S_t^F S_s^F \Bu = S_{t+s}^F \Bu$ for all $\Bu \in \mathscr{D}^F$, $t,s \geq 0$.
    \item \textbf{Lipschitz continuity}: There exists a constant $L_F > 0$ such that for all $\Bu, \boldsymbol{w} \in \mathscr{D}^F$ and $t,s \geq 0$,
    \[
    \|S_t^F \Bu - S_s^F \boldsymbol{w}\|_{L^1} \leq L_F \left(|t-s| + \|\Bu - \boldsymbol{w}\|_{L^1}\right).
    \]
    \item \textbf{Riemann consistency}: For every piecewise constant initial data $\Bu \in \mathscr{D}^F$, the semigroup solution $S_t^F \Bu$ coincides with the gluing of standard Riemann problem solutions for sufficiently small $t > 0$.
\end{itemize}
Here $\mathscr{D}^F \subset \mathbf{L}^1(\mathbb{R}; \mathbb{R}^2)$ denotes the domain containing integrable functions with sufficiently small total variation.
\end{Definition}

For all $\Bu \in \mathscr{D}^F$, the map $t \mapsto S_t^F \Bu$ provides the unique entropic solution to the Cauchy problem \cite{BA1995}. The stability theory developed in \cite{BSCR} ensures that small perturbations in the flux function lead to proportionally small changes in the semigroup, which provides the theoretical foundation for the convergence analysis of our flux reconstruction method.

The characteristic structure of the system is encoded in the Jacobian matrix:
\begin{equation*}
DF(\Bu) = \begin{pmatrix}
0 & f_1'(v) \\
f_2'(u) & 0
\end{pmatrix},
\end{equation*}
which possesses eigenvalues $\lambda_{1,2}(\Bu) = \mp \sqrt{f_1'(v)f_2'(u)}$ and corresponding eigenvectors. This eigenstructure governs the propagation of information through the system and forms the mathematical foundation for our wave front tracking approach to flux inversion.

The fundamental building block of our analysis is the Riemann problem with piecewise constant initial data:
\begin{equation}\label{eq:initial}
  \Bu(0, x)= \begin{cases}\Bu_l, & x>0 ; \\ \Bu_r, & x<0 .\end{cases}
\end{equation}
Here $\Bu_l=\left(\begin{array}{c}u_l \\ v_l\end{array}\right)$ and $\Bu_r=\left(\begin{array}{c}u_r \\ v_r\end{array}\right)$ represent constant states separated by an initial discontinuity. As concisely summarized by Courant and Friedrichs \cite{CM2014}, Riemann solutions exhibit characteristic wave patterns, ``either the initial discontinuity is resolved immediately and the disturbance, while propagated, becomes continuous, or the initial discontinuity is propagated through one or two shock fronts, advancing not at sonic but at supersonic speed relative to the medium ahead of them." This fundamental dichotomy between continuous and discontinuous wave propagation, governed by entropy considerations, provides the key insight that enables our flux reconstruction methodology.

To be more precise, there are four fundamental solution patterns for the Riemann problem, which connect constant states through combinations of shocks and rarefaction waves. We now describe the wave phenomena that constitute these solution patterns.

Shock waves are piecewise constant discontinuous solutions that propagate initial discontinuities while satisfying the Rankine-Hugoniot jump conditions and entropy inequalities across the shock front. For the coupled system \eqref{Main:eqf}, there exist two families of shock waves associated with the characteristic speeds $\lambda_1(\Bu) = -\sqrt{f_1'(v)f_2'(u)}$ and $\lambda_2(\Bu) = \sqrt{f_1'(v)f_2'(u)}$. We focus on the second family (front shocks), as the analysis for the first family follows by symmetry.

A second-family shock is a solution consisting of two constant states $\Bu_l = (u_l, v_l)$ and $\Bu_r = (u_r, v_r)$ separated by a ray $x = st$, where the shock speed $s$ is determined by the Rankine-Hugoniot conditions:
\begin{equation}\label{eq:RHC}
\begin{cases}
f_1(v_r) - f_1(v_l) = s(u_r - u_l), \\
f_2(u_r) - f_2(u_l) = s(v_r - v_l).
\end{cases}
\end{equation}

Eliminating the shock speed $s$ yields the Hugoniot locus:
\begin{equation*}
\frac{f_1(v_r) - f_1(v_l)}{u_r - u_l} = \frac{f_2(u_r) - f_2(u_l)}{v_r - v_l}.
\end{equation*}

\begin{figure}[htbp]
  \centering
  \includegraphics[scale=0.47]{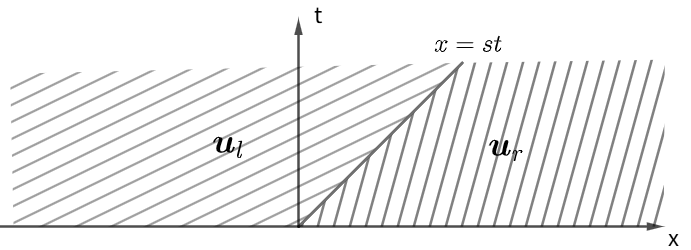}
  \caption{Second-family shock wave (front shock)} \label{fg:01}
\end{figure}

For a genuinely nonlinear second family, the shock satisfies the Lax entropy condition:
\begin{equation*}
\lambda_2(\Bu_l) < s < \lambda_2(\Bu_r),
\end{equation*}
which ensures the admissibility of the discontinuous solution. This entropy condition guarantees that characteristics converge into the shock front, providing a physical selection criterion for physically meaningful solutions. Figure \ref{fg:01} provides a visual illustration of a second-family shock wave.

In contrast to shock waves, rarefaction waves are continuous self-similar solutions that smoothly connect constant states, taking the form $\Bu(t,x) = \Bu(\xi)$ with $\xi = x/t$. For the coupled system \eqref{Main:eqf}, there are two families of rarefaction waves: the first family (back rarefactions) associated with the negative eigenvalue $\lambda_1(\Bu) = -\sqrt{f_1'(v)f_2'(u)}$, and the second family (front rarefactions) associated with the positive eigenvalue $\lambda_2(\Bu) = \sqrt{f_1'(v)f_2'(u)}$. We focus on the second family, as the first family follows by symmetry.

For a second-family rarefaction wave, the self-similar ansatz $\Bu(t,x) = \Bu(\xi)$ reduces system \eqref{Main:eqf} to the ordinary differential equation:
\begin{equation*}
(-\xi I + A(\Bu))\frac{d\Bu}{d\xi} = 0,
\end{equation*}
where $A(\Bu) = \begin{pmatrix} 0 & f_1'(v) \\ f_2'(u) & 0 \end{pmatrix}$ is the Jacobian matrix.

This implies that $\Bu(\xi)$ follows the integral curve of the second characteristic field:
\begin{equation*}
\frac{du}{dv} = \frac{f_1'(v)}{\lambda_2(\Bu)} = \sqrt{\frac{f_1'(v)}{f_2'(u)}}.
\end{equation*}

The rarefaction wave connecting states $\Bu_l = (u_l, v_l)$ and $\Bu_r = (u_r, v_r)$ satisfies the integral relation:
\begin{equation}\label{eq:RW}
u_r - u_l = r_2(v_r; \Bu_l) := \int_{v_l}^{v_r} \sqrt{\frac{f_1'(w)}{f_2'(u(w))}} dw, \quad v_r < v_l,
\end{equation}
where $u(w)$ follows the integral curve from $\Bu_l$.

\begin{figure}[htbp]
  \centering
  \includegraphics[scale=0.45]{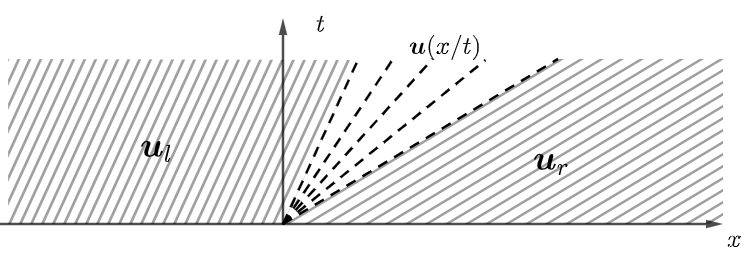}
  \caption{Second-family rarefaction wave (front rarefaction)} 
  \label{fg:02}
\end{figure}

The monotonicity condition $v_r < v_l$ ensures that the characteristic speed $\lambda_2(\Bu(\xi))$ increases with $\xi$, satisfying the genuine nonlinearity requirement. The rarefaction fan spans the interval $\xi \in [\lambda_2(\Bu_l), \lambda_2(\Bu_r)]$, with the solution constant along each ray $\xi = \text{constant}$. Figure \ref{fg:02} illustrates the structure of a second-family rarefaction wave.

The centered rarefaction region is characterized by rays emanating from the origin. The left boundary $\xi = \lambda_2(\Bu_l)$ is determined by the left state $\Bu_l$, while the right boundary $\xi = \lambda_2(\Bu_r)$ matches the characteristic speed of the right state $\Bu_r$. Within this fan, the solution varies smoothly along the integral curve of the second characteristic field.

\begin{figure}[htbp]
  \centering
  \includegraphics[scale=0.3]{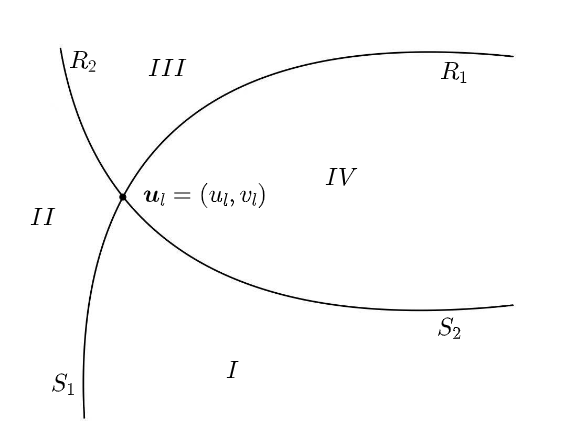}
  \caption{Wave pattern classification: Regions I-IV } \label{fg:03}
\end{figure}

In the $(v, u)$-plane, for a given point $\Bu_l$, the equation \eqref{eq:RHC} and its counterpart for front shocks define two curves $S_1$ and $S_2$ emanating from the point $\Bu_l$; similarly \eqref{eq:RW} and its counterpart for front rarefaction waves also define two curves $R_1$ and $R_2$ emanating from the point $\Bu_l$. We now study the four fundamental solution patterns for the Riemann problem. For any given point $\Bu_l$ in the $(u, v)$-plane, for $i=1,2$, the curves $S_i$ and $R_i$ join at $\Bu_l$. The two intersecting curves divide the $(u, v)$-plane into four parts. We list them as I, II, III and IV (see Fig. \ref{fg:03}):

\begin{align*}
    I:& \text{bordering on } S_1,S_2, \quad &II: \text{bordering on } S_1,R_2,\\
    III:& \text{bordering on } R_1,R_2, \quad &IV: \text{bordering on } R_1,S_2.
\end{align*}

For $\Bu_r$ sufficiently close to $\Bu_l$, the structure of the general solution to the Riemann problem \eqref{Main:eqf}-\eqref{eq:initial} is now determined by the location of the state $\Bu_r$ with respect to the curves $R_i,S_i$ for $i=1,2$. We will analyze these four cases in detail in Section \ref{Sec:proof} and use this feature to give the reconstruction algorithm of the stream functions $f_1$ and $f_2$.

\subsection{Main Results}
Building upon the mathematical framework established in the previous section, we now present the main theoretical contributions of this work. Our central result establishes that under reasonable regularity assumptions, it is possible to reconstruct the unknown flux functions from single-time observations of Riemann solutions. This represents a significant advancement in the theory of inverse problems for hyperbolic conservation laws, as it circumvents the need for dense temporal data that has limited previous approaches.

We begin by formalizing the concept of observability that underpins our reconstruction methodology:
\begin{Definition}\label{Def1}
   A function $z:[0, \infty) \times \mathbb{R} \rightarrow \mathbb{R}$ is said to be observable at fixed $T>0$ if we know its values $z(x,T)$ for (almost) every $x \in \mathbb{R}$.
\end{Definition}

This observability condition captures the practical scenario where we have access to solution profiles at a specific observation time, but not necessarily throughout the temporal evolution. Under this condition, our main reconstruction theorem demonstrates the feasibility of flux identification:

\begin{thm}\label{main:thm-f}
  Let $T>0$ and let $\Omega \subset \mathbb{R}^2$ be an open domain containing the rectangle $\left[u_*, u^*\right] \times \left[v_*, v^*\right]$, where $u_*<u^*$ and $v_*<v^*$. Let $c_1, c_2 \in \mathbb{R}$ be fixed. Assume that the flux functions $f_1, f_2: \mathbb{R} \rightarrow \mathbb{R}$ are $C^{1,1}$ with finite numbers of inflection points, that $f_1(v_*)=c_1$, $f_2(u_*)=c_2$, and that the coupled flux $F(\Bu) = (f_1(v), f_2(u))^T$ generates a Standard Riemann Semigroup on $\Omega$. Furthermore, assume that the solution to any Riemann problem for \eqref{Main:eqf} at time $T$ is observable, in the sense of definition \ref{Def1}. 
  
  Then, for all $m \in \mathbb{N}$, setting $\delta:=2^{-m}(u^*-u_*)$, $\eta:=2^{-m}(v^*-v_*)$ and $u_\alpha=u_*+\alpha \delta$, $v_\alpha=v_*+\alpha \eta$ for $\alpha=0, \ldots, 2^m$, there exist two piecewise quadratic $C^1$ functions $f_{1,m}:\left[v_*, v^*\right] \rightarrow \mathbb{R}$, $f_{2,m}:\left[u_*, u^*\right] \rightarrow \mathbb{R}$ such that $f_{1,m}(v_\alpha)=f_{1}(v_\alpha)$ and $f_{2,m}(u_\alpha)=f_{2}(u_\alpha)$ for all $\alpha$, and the following error estimates hold:
\begin{equation}\label{eq:result-f}
  \begin{split}
    \|f_{1} - f_{1,m}\|_{L^\infty([v_*,v^*])} &\leqslant L_1\eta^2,\\
    \|f_{2} - f_{2,m}\|_{L^\infty([u_*,u^*])} &\leqslant L_2\delta^2,\\
    \|f_{1}' - f_{1,m}'\|_{L^\infty([v_*,v^*])} &\leqslant 3L_1\eta,\\
    \|f_{2}' - f_{2,m}'\|_{L^\infty([u_*,u^*])} &\leqslant 3L_2\delta,
  \end{split}
\end{equation}
where $L_1$ and $L_2$ are the Lipschitz constants of $f_1'$ and $f_2'$ on the respective intervals.

The practical significance of this mathematical reconstruction is validated by the following stability result: if $\hat{\Bu}$ is a $\mathbf{BV}$ function with values in $\left[u_*, u^*\right] \times \left[v_*, v^*\right]$, and we denote by $\Bu^m$ (resp. $\Bu_{\mathrm{obs}}$) the solution to the Cauchy problem with flux $(f_{1,m}, f_{2,m})$ (resp. $(f_1, f_2)$) with initial data $\hat{\Bu}$, then
\begin{equation}\label{eq:stability}
    \left\|S_{T}^{F_m} \Bu^m-S_{T}^{F}\Bu_{\mathrm{obs}}\right\|_{\mathbf{L}^1} \leqslant C T (L_1\eta + L_2\delta),
\end{equation}
where $F_m=(f_{1,m},f_{2,m})$, $F=(f_1,f_2)$, and $C$ is a constant that does not depend on $\delta$ and $\eta$.
\end{thm}

The power of Theorem \ref{main:thm-f} lies in its demonstration that accurate flux reconstruction is achievable from limited observational data. The error estimates provide explicit guidance on how the discretization parameters affect reconstruction accuracy, while the stability result ensures that small errors in flux approximation translate to proportionally small errors in solution prediction.

\begin{remark}
\textbf{(Invariance under Reference Value Shifts).} 
Even when the reference values $c_1 = f_1(v_*)$ and $c_2 = f_2(u_*)$ are unknown, our reconstruction algorithm remains effective with only minor modifications to the error estimates. In this scenario, the reconstructed fluxes $f_{1,m}$ and $f_{2,m}$ will differ from the true fluxes by at most an additive constant:
\[
\|f_{1} - f_{1,m}\|_{L^\infty([v_*,v^*])} \leq L_1\eta^2 + C_1, \quad \|f_{2} - f_{2,m}\|_{L^\infty([u_*,u^*])} \leq L_2\delta^2 + C_2,
\]
where $C_1, C_2$ are constants representing the uncertainty in the reference values. Crucially, the derivative estimates remain unaffected:
\[
\|f_{1}' - f_{1,m}'\|_{L^\infty([v_*,v^*])} \leq 3L_1\eta, \quad \|f_{2}' - f_{2,m}'\|_{L^\infty([u_*,u^*])} \leq 3L_2\delta,
\]
and the stability estimate \eqref{eq:stability} for the solution operators remains valid with the same constant $C$. This invariance property stems from the fundamental structure of hyperbolic conservation laws: the characteristic speeds, shock speeds, and wave interactions depend only on flux derivatives rather than absolute flux values. Consequently, our wave front tracking methodology accurately recovers the essential dynamical information even in the presence of reference value uncertainty.
\end{remark}

An important refinement of these results emerges when the flux functions possess higher regularity, as captured in the following remark:
\begin{remark}
    \textbf{(Enhanced Convergence under $C^3$ Regularity).} 
    If the flux functions possess higher regularity, specifically if $f_1 \in C^3[v_*, v^*]$ and $f_2 \in C^3[u_*, u^*]$ with bounded third derivatives, then the approximation quality of our piecewise quadratic $C^1$ interpolants improves significantly. In this case, the convergence rates enhance to:
    \begin{align*}
        \|f_1 - f_{1,m}\|_{L^\infty([v_*,v^*])} &= O(\eta^3), \quad \|f_1' - f_{1,m}'\|_{L^\infty([v_*,v^*])} = O(\eta^2),\\
        \|f_2 - f_{2,m}\|_{L^\infty([u_*,u^*])} &= O(\delta^3), \quad \|f_2' - f_{2,m}'\|_{L^\infty([u_*,u^*])} = O(\delta^2),
    \end{align*}
    and consequently, the solution error estimate improves to:
    \begin{equation*}
        \left\|S_{T}^{F_m} \Bu^m-S_{T}^{F}\Bu_{\mathrm{obs}}\right\|_{\mathbf{L}^1} \leqslant C T (L_1\eta^2 + L_2\delta^2).
    \end{equation*}
    This remarkable improvement demonstrates that our reconstruction method benefits substantially from higher regularity of the unknown flux functions, achieving rapid convergence in practical applications where the physical fluxes are smooth.
\end{remark}
Together, these results establish a comprehensive framework for flux identification in coupled hyperbolic systems. The combination of rigorous error estimates, stability guarantees, and enhanced convergence under higher regularity provides both theoretical foundation and practical guidance for implementing our wave front tracking methodology in various applied contexts.

\section{Proofs of Main Results}\label{Sec:proof}

\subsection{Wave Front Tracking Analysis}

We now proceed to construct the approximated flux functions $f_{1,m}$ and $f_{2,m}$ that satisfy the error estimates (\ref{eq:result-f}). 

Let $T > 0$, $u_*, u^* \in \mathbb{R}$ with $u_* < u^*$, $v_*, v^* \in \mathbb{R}$ with $v_* < v^*$, and $c_1, c_2 \in \mathbb{R}$ be fixed parameters. Fix $m \in \mathbb{N}$ and define the mesh sizes $\delta = 2^{-m}|u^* - u_*|$, $\eta = 2^{-m}|v^* - v_*|$, and the grid points:
\begin{equation*}
    u_\alpha = u_* + \alpha \delta, \, v_\alpha = v_* +\alpha\eta \quad \text{for } \alpha = 0, \ldots, 2^m.
\end{equation*}

We initialize the approximations by setting $f_{1,m}(v_*) = c_1$ and $f_{2,m}(u_*) = c_2$.

Since our only available information comes from observing solutions to the coupled system (\ref{Main:eqf}) at the fixed time $T > 0$, we must carefully design the initial data to extract information about the unknown flux functions $f_1$ and $f_2$. To this end, we consider the following family of Riemann initial data:

\begin{equation*}
    \Bu_0^h(x) = \begin{cases}
        \Bu_h, & x < 0, \\
        \Bu_{h+1}, & x > 0,
    \end{cases} \quad h = 0, \ldots, 2^m - 1,
\end{equation*}
where $\Bu_h = \begin{pmatrix} u_h \\ v_h \end{pmatrix}$. The reconstruction strategy is iterative: we use the solution corresponding to $\Bu_0^h$ to determine $f_{1,m}$ at $v_{h+1}$ and $f_{2,m}$ at $u_{h+1}$.

The solutions to these Riemann problems exhibit four fundamental wave patterns, each connecting constant states through combinations of shocks and rarefaction waves in the two characteristic families. We now analyze each case in detail.

\medskip

\begin{figure}[htbp]
    \centering
     \begin{subfigure}[b]{0.3\textwidth}  
        \centering
        \includegraphics[width=\linewidth]{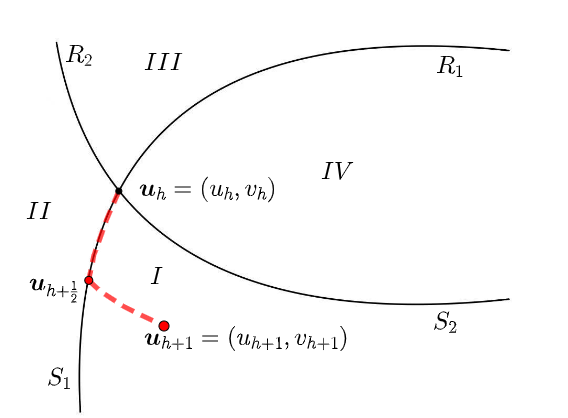}
        \caption{Phase plane: states configuration}
    \end{subfigure}
     \hspace{10mm}  
    \begin{subfigure}[b]{0.45\textwidth} 
        \centering
        \includegraphics[width=\linewidth]{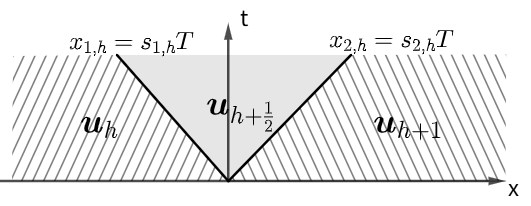}
       \caption{Wave structure at time $T$}
    \end{subfigure}
    \caption{Two shock waves connecting states $(u_h,v_h)$ and $(u_{h+1},v_{h+1})$ through intermediate state $(u_{h+\frac{1}{2}},v_{h+\frac{1}{2}})$ for system \eqref{Main:eqf}.}
    \label{fig:ss}
\end{figure}

\textbf{Case 1 (Two Shocks).} 
For the coupled system \eqref{Main:eqf}, shock waves are piecewise constant discontinuous solutions propagating initial discontinuities while satisfying the Rankine-Hugoniot conditions. We consider the case where the Riemann solution consists of two shock waves: a first-family shock (associated with $\lambda_1(\Bu) = -\sqrt{f_1'(v)f_2'(u)}$) followed by a second-family shock (associated with $\lambda_2(\Bu) = \sqrt{f_1'(v)f_2'(u)}$), as illustrated in Figure \ref{fig:ss}.

For a fixed $h \geq 0$, let $\tilde{\Bu}(\cdot) = \begin{pmatrix} \tilde{u}(\cdot) \\ \tilde{v}(\cdot) \end{pmatrix} = \Bu_{\text{obs}}(T, \cdot)$ be the observed solution at time $T$, consisting of two shock waves connecting the states $(u_h, v_h)$ and $(u_{h+1}, v_{h+1})$. Let $x_{1,h}$ and $x_{2,h}$ denote the positions of the first and second shocks, respectively, with corresponding propagation speeds $s_{1,h} = x_{1,h}/T$ and $s_{2,h} = x_{2,h}/T$. Let $(u_{h+\frac{1}{2}}, v_{h+\frac{1}{2}})$ be the intermediate constant state between the two shocks.

The Rankine-Hugoniot conditions across both shocks yield:
\begin{equation*}
    \begin{cases}
        f_1(v_{h+\frac{1}{2}}) - f_1(v_h) = s_{1,h}(u_{h+\frac{1}{2}} - u_h), \\
        f_2(u_{h+\frac{1}{2}}) - f_2(u_h) = s_{1,h}(v_{h+\frac{1}{2}} - v_h), \\
        f_1(v_{h+1}) - f_1(v_{h+\frac{1}{2}}) = s_{2,h}(u_{h+1} - u_{h+\frac{1}{2}}), \\
        f_2(u_{h+1}) - f_2(u_{h+\frac{1}{2}}) = s_{2,h}(v_{h+1} - v_{h+\frac{1}{2}}).
    \end{cases}
\end{equation*}
It yeids 
 \begin{equation*}
    \left\{\begin{array}{l}
        f_{1}(v_{h+1})  = f_{1}(v_{h}) + s_{1,h}(u_{h+\frac{1}{2}}-u_{h}) + s_{2,h}(u_{h+1}-u_{h+\frac{1}{2}}),\\
        f_{2}(u_{h+1})  = f_{2}(u_{h}) + s_{1,h}(v_{h+\frac{1}{2}}-v_{h}) + s_{2,h}(v_{h+1}-v_{h+\frac{1}{2}}).
      \end{array}\right.
  \end{equation*} 
  Therefore, if $f_{1,m}(v_0),\ldots , f_{1,m}(v_h)$ and $ f_{2,m}(u_0),\ldots , f_{2,m}(u_h)$ are given so that $f_{1,m}(v_{\alpha}) = f_{1}(v_{\alpha})$ and $f_{2,m}(u_{\alpha}) = f_{2}(u_{\alpha})$, we can define
  \begin{equation}
    \left\{\begin{array}{l}
      f_{1,m}(v_{h+1})  := f_{1,m}(v_{h}) + s_{1,h}(u_{h+\frac{1}{2}}-u_{h}) + s_{2,h}(u_{h+1}-u_{h+\frac{1}{2}})= f_{1}(v_{h+1}),\\
      f_{2,m}(u_{h+1})  := f_{2,m}(u_{h}) + s_{1,h}(v_{h+\frac{1}{2}}-v_{h}) + s_{2,h}(v_{h+1}-v_{h+\frac{1}{2}}) = f_{2}(u_{h+1}).
      \end{array}\right.
  \end{equation} 

\medskip

\textbf{Case 2 (Shock + Rarefaction).} 

\begin{figure}[htbp]
    \centering
     \begin{subfigure}[b]{0.3\textwidth}  
        \centering
        \includegraphics[width=\linewidth]{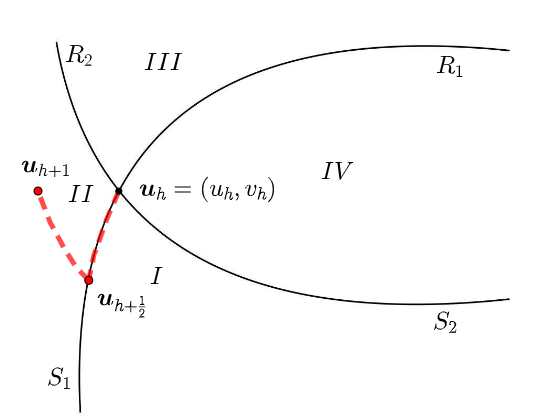}
        \caption{Phase plane: states configuration}
    \end{subfigure}
     \hspace{10mm}  
    \begin{subfigure}[b]{0.45\textwidth} 
        \centering
        \includegraphics[width=\linewidth]{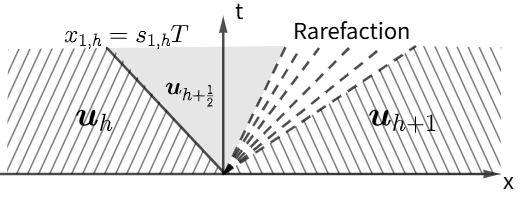}
       \caption{Wave structure at time $T$}
    \end{subfigure}
    \caption{First-family shock wave followed by second-family rarefaction wave for system \eqref{Main:eqf}.}
    \label{fig:sr}
\end{figure}

For the coupled system \eqref{Main:eqf}, we now consider the case where the Riemann solution consists of a first-family shock wave followed by a second-family rarefaction wave. This wave pattern connects the left state $(u_h, v_h)$ to the right state $(u_{h+1}, v_{h+1})$ through an intermediate state $(u_{h+\frac{1}{2}}, v_{h+\frac{1}{2}})$, as illustrated in Figure \ref{fig:sr}.

To incorporate rarefaction waves within our reconstruction framework, we introduce the concept of \emph{equivalent shock}. The fundamental idea is to replace a continuous rarefaction fan with a single discontinuity positioned to preserve essential integral quantities derived from the conservation laws. For a second-family rarefaction wave in the $u$-component connecting states $u_{h+\frac{1}{2}}$ and $u_{h+1}$ over the spatial interval $I_{2,h} = [x_{2,h}, x_{2,h+1}]$, we define the equivalent shock position as the mass centroid:
\begin{equation}\label{eq:es-u}
    \xi_{u,h} := \frac{\int_{u_{h+\frac{1}{2}}}^{u_{h+1}} x_1(u) du}{u_{h+1} - u_{h+\frac{1}{2}}},
\end{equation}
where $x_1(u)$ is the inverse function of $\tilde{u}(x)$ on $I_{2,h}$. This positioning ensures conservation of the average spatial coordinate for each level set of $u$. The corresponding equivalent shock speed is $s_{u,h} = \xi_{u,h}/T$. Similarly, for the $v$-component within the same rarefaction fan connecting $v_{h+\frac{1}{2}}$ and $v_{h+1}$, we define:
\begin{equation}
    \xi_{v,h} := \frac{\int_{v_{h+\frac{1}{2}}}^{v_{h+1}} x_2(v) dv}{v_{h+1} - v_{h+\frac{1}{2}}},
\end{equation}
with equivalent shock speed $s_{v,h} = \xi_{v,h}/T$.

The mathematical justification for the equivalent shock approximation stems from the intrinsic structure of centered rarefaction waves. For a self-similar solution of the form $\Bu(x/t)$, the conservation laws impose specific integral constraints. Consider the first conservation law $\partial_t u + \partial_x f_1(v) = 0$ over the rarefaction region $I_{2,h}$. At fixed time $T$, integrating over the spatial domain yields:
\[
\int_{x_{2,h}}^{x_{2,h+1}} [\partial_t u + \partial_x f_1(v)] dx = 0.
\]
Using the self-similar form with $\xi = x/T$, we have $\partial_t u = -\frac{x}{T^2}u'(\xi)$ and $\partial_x f_1(v) = \frac{1}{T}f_1'(v)v'(\xi)$, which gives:
\[
-\frac{1}{T^2} \int_{x_{2,h}}^{x_{2,h+1}} x u'(\xi) dx + \frac{1}{T}[f_1(v(x_{2,h+1})) - f_1(v(x_{2,h}))] = 0.
\]
Changing variables to $\xi = x/T$ with $dx = T d\xi$, and noting that $u'(\xi)d\xi = du$, we obtain:
\[
-\frac{1}{T} \int_{\xi_{2,h}}^{\xi_{2,h+1}} \xi du + [f_1(v_{h+1}) - f_1(v_{h+\frac{1}{2}})] = 0.
\]
Since $\xi = x/T$ and $x = x_1(u)$ within the rarefaction, this becomes:
\begin{equation}
    \int_{u_{h+\frac{1}{2}}}^{u_{h+1}} \frac{x_1(u)}{T} du = f_1(v_{h+1}) - f_1(v_{h+\frac{1}{2}}).
\end{equation}

Now consider the second conservation law $\partial_t v + \partial_x f_2(u) = 0$ over the same region:
\[
\int_{x_{2,h}}^{x_{2,h+1}} [\partial_t v + \partial_x f_2(u)] dx = 0.
\]
Following a similar derivation with $v'(\xi)d\xi = dv$, we obtain:
\begin{equation}\label{eq:sss-v}
    \int_{v_{h+\frac{1}{2}}}^{v_{h+1}} \frac{x_2(v)}{T} dv = f_2(u_{h+1}) - f_2(u_{h+\frac{1}{2}}).
\end{equation}

These integral relations provide the rigorous foundation for the equivalent shock approximation in general coupled systems.

For fixed $h \geq 0$, let $\tilde{\Bu}(\cdot) = \Bu_{\text{obs}}(T, \cdot)$ be the observed solution consisting of a first-family shock followed by a second-family rarefaction connecting $(u_h, v_h)$ to $(u_{h+1}, v_{h+1})$. Let $x_{1,h}$ be the shock position with speed $s_{1,h} = x_{1,h}/T$, and let $(u_{h+\frac{1}{2}}, v_{h+\frac{1}{2}})$ be the intermediate state.

The Rankine-Hugoniot conditions across the first-family shock give:
\begin{equation}\label{eq:SR1}
    \begin{cases}
        f_1(v_{h+\frac{1}{2}}) - f_1(v_h) = s_{1,h}(u_{h+\frac{1}{2}} - u_h), \\
        f_2(u_{h+\frac{1}{2}}) - f_2(u_h) = s_{1,h}(v_{h+\frac{1}{2}} - v_h).
    \end{cases}
\end{equation}

For the second-family rarefaction region $I_{2,h} = [x_{2,h}, x_{2,h+1}]$, applying the equivalent shock approximation \eqref{eq:es-u}-\eqref{eq:sss-v} yields:
\begin{equation}\label{eq:SR2}
    \begin{cases}
        f_1(v_{h+1}) - f_1(v_{h+\frac{1}{2}}) = s_{u,h}(u_{h+1} - u_{h+\frac{1}{2}}), \\
        f_2(u_{h+1}) - f_2(u_{h+\frac{1}{2}}) = s_{v,h}(v_{h+1} - v_{h+\frac{1}{2}}).
    \end{cases}
\end{equation}

From equations (\ref{eq:SR1}) and (\ref{eq:SR2}), we can directly reconstruct the flux functions at the intermediate and right states. The key observation is that the intermediate state values $u_{h+\frac{1}{2}}$ and $v_{h+\frac{1}{2}}$ are known from the observed solution, allowing us to eliminate the intermediate flux values $f_1(v_{h+\frac{1}{2}})$ and $f_2(u_{h+\frac{1}{2}})$. Combine \eqref{eq:SR1} with \eqref{eq:SR2}, we can get
\begin{equation*}
    \begin{cases}
        f_1(v_{h+1}) = f_1(v_h) + s_{1,h}(u_{h+\frac{1}{2}} - u_h) + s_{u,h}(u_{h+1} - u_{h+\frac{1}{2}}),\\
        f_2(u_{h+1}) = f_2(u_h) + s_{1,h}(v_{h+\frac{1}{2}} - v_h) + s_{v,h}(v_{h+1} - v_{h+\frac{1}{2}}).
    \end{cases}
\end{equation*}

Therefore, if $f_{1,m}(v_0), \ldots, f_{1,m}(v_h)$ and $f_{2,m}(u_0), \ldots, f_{2,m}(u_h)$ are known with $f_{1,m}(v_\beta) = f_1(v_\beta)$ and $f_{2,m}(u_\alpha) = f_2(u_\alpha)$, we extend the approximations by:
\begin{equation}
    \begin{cases}
        f_{1,m}(v_{h+1}) = f_{1,m}(v_h) + s_{1,h}(u_{h+\frac{1}{2}} - u_h) + s_{u,h}(u_{h+1} - u_{h+\frac{1}{2}}) = f_1(v_{h+1}), \\
        f_{2,m}(u_{h+1}) = f_{2,m}(u_h) + s_{1,h}(v_{h+\frac{1}{2}} - v_h) + s_{v,h}(v_{h+1} - v_{h+\frac{1}{2}}) = f_2(u_{h+1}).
    \end{cases}
\end{equation}

This reconstruction strategy successfully eliminates the intermediate flux values $f_1(v_{h+\frac{1}{2}})$ and $f_2(u_{h+\frac{1}{2}})$ by leveraging the directly observable intermediate state $(u_{h+\frac{1}{2}}, v_{h+\frac{1}{2}})$ and the measured wave speeds $s_{1,h}$, $s_{u,h}$, and $s_{v,h}$.

\medskip

\textbf{Case 3 (Rarefaction + Shock).} 

\begin{figure}[htbp]
    \centering
    \begin{subfigure}[b]{0.3\textwidth}
        \centering
        \includegraphics[width=\linewidth]{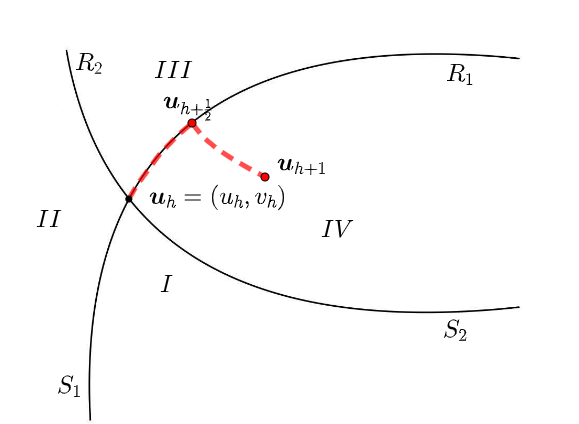}
        \caption{Phase plane: states configuration}
    \end{subfigure}
    \hspace{10mm} 
    \begin{subfigure}[b]{0.45\textwidth}
        \centering
        \includegraphics[width=\linewidth]{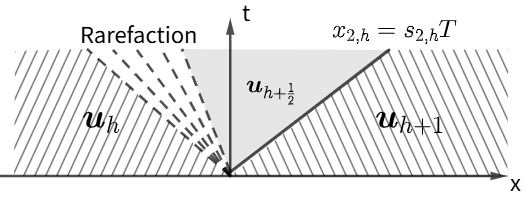}
        \caption{Wave structure at time $T$}
    \end{subfigure}
    \caption{First-family rarefaction followed by second-family shock wave for system \eqref{Main:eqf}.}
    \label{fig:rs}
\end{figure}

This case exhibits symmetry with Case 2, featuring a first-family rarefaction wave followed by a second-family shock wave. For fixed $h \geq 0$, let $\tilde{\Bu}(\cdot) = \Bu_{\text{obs}}(T, \cdot)$ be the observed solution connecting states $(u_h, v_h)$ and $(u_{h+1}, v_{h+1})$ through an intermediate state $(u_{h+\frac{1}{2}}, v_{h+\frac{1}{2}})$. Let $x_{2,h}$ be the shock position with speed $s_{2,h} = x_{2,h}/T$, as illustrated in Figure \ref{fig:rs}.

The Rankine-Hugoniot conditions across the second-family shock give:
\begin{equation}\label{eq:RS1}
    \begin{cases}
        f_1(v_{h+1}) - f_1(v_{h+\frac{1}{2}}) = s_{2,h}(u_{h+1} - u_{h+\frac{1}{2}}), \\
        f_2(u_{h+1}) - f_2(u_{h+\frac{1}{2}}) = s_{2,h}(v_{h+1} - v_{h+\frac{1}{2}}).
    \end{cases}
\end{equation}

For the first-family rarefaction region $I_{1,h} = [x_{1,h}, x_{1,h+1}]$ connecting $(u_h, v_h)$ to $(u_{h+\frac{1}{2}}, v_{h+\frac{1}{2}})$, we apply the equivalent shock approximation. Define the equivalent shock positions:
\begin{equation*}
    \xi_{u,h} := \frac{\int_{u_h}^{u_{h+\frac{1}{2}}} x_1(u) du}{u_{h+\frac{1}{2}} - u_h}, \quad
    \xi_{v,h} := \frac{\int_{v_h}^{v_{h+\frac{1}{2}}} x_2(v) dv}{v_{h+\frac{1}{2}} - v_h},
\end{equation*}
with equivalent shock speeds $s_{u,h} = \xi_{u,h}/T$ and $s_{v,h} = \xi_{v,h}/T$.

From the integral relations for centered rarefaction waves in the general system \eqref{Main:eqf}, we obtain:
\begin{equation}\label{eq:RS2}
    \begin{cases}
        f_1(v_{h+\frac{1}{2}}) - f_1(v_h) = s_{u,h}(u_{h+\frac{1}{2}} - u_h), \\
        f_2(u_{h+\frac{1}{2}}) - f_2(u_h) = s_{v,h}(v_{h+\frac{1}{2}} - v_h).
    \end{cases}
\end{equation}

Using the observable intermediate state $(u_{h+\frac{1}{2}}, v_{h+\frac{1}{2}})$ and the measured wave speeds, we can directly reconstruct the flux functions. Substituting equation (\ref{eq:RS2}) into equation (\ref{eq:RS1}) gives the reconstruction for the right state:
\begin{equation*}
    \begin{cases}
        f_1(v_{h+1}) = f_1(v_h) + s_{u,h}(u_{h+\frac{1}{2}} - u_h) + s_{2,h}(u_{h+1} - u_{h+\frac{1}{2}}), \\
        f_2(u_{h+1}) = f_2(u_h) + s_{v,h}(v_{h+\frac{1}{2}} - v_h) + s_{2,h}(v_{h+1} - v_{h+\frac{1}{2}}).
    \end{cases}
\end{equation*}

Therefore, if $f_{1,m}(v_0), \ldots, f_{1,m}(v_h)$ and $f_{2,m}(u_0), \ldots, f_{2,m}(u_h)$ are known with $f_{1,m}(v_\beta) = f_1(v_\beta)$ and $f_{2,m}(u_\alpha) = f_2(u_\alpha)$, we extend the approximations by:
\begin{equation}
    \begin{cases}
        f_{1,m}(v_{h+1}) = f_{1,m}(v_h) + s_{u,h}(u_{h+\frac{1}{2}} - u_h) + s_{2,h}(u_{h+1} - u_{h+\frac{1}{2}}) = f_1(v_{h+1}), \\
        f_{2,m}(u_{h+1}) = f_{2,m}(u_h) + s_{v,h}(v_{h+\frac{1}{2}} - v_h) + s_{2,h}(v_{h+1} - v_{h+\frac{1}{2}}) = f_2(u_{h+1}).
    \end{cases}
\end{equation}

This reconstruction successfully utilizes the observable intermediate state and wave speeds to determine both flux functions $f_1$ and $f_2$ at the grid points, maintaining the symmetric relationship with Case 2 while accounting for the reversed wave order.

\medskip

\textbf{Case 4 (Two Rarefactions).} 

\begin{figure}[htbp]
    \centering
    \begin{subfigure}[b]{0.3\textwidth}
        \centering
        \includegraphics[width=\linewidth]{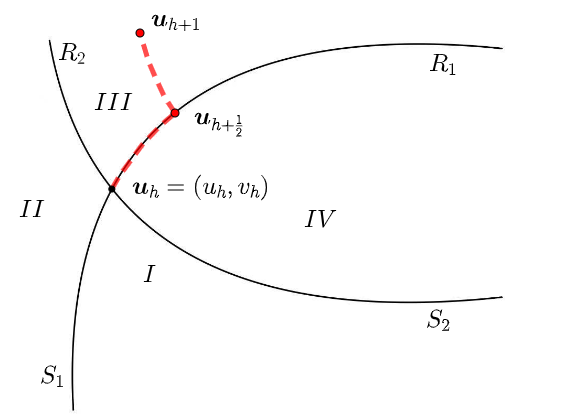}
        \caption{Phase plane: states configuration}
    \end{subfigure}
    \hspace{10mm} 
    \begin{subfigure}[b]{0.45\textwidth}
        \centering
        \includegraphics[width=\linewidth]{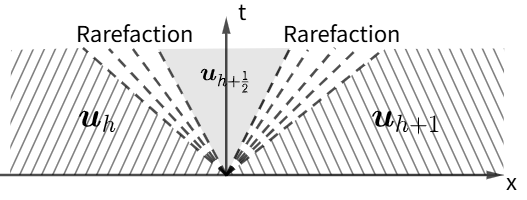}
        \caption{Wave structure at time $T$}
    \end{subfigure}
    \caption{First-family rarefaction followed by second-family rarefaction for system \eqref{Main:eqf}.}
    \label{fig:rr}
\end{figure}

For fixed $h \geq 0$, let $\tilde{\Bu}(\cdot) = \Bu_{\text{obs}}(T, \cdot)$ be the observed solution consisting of two rarefaction waves: a first-family rarefaction followed by a second-family rarefaction. The solution connects the left state $(u_h, v_h)$ to the right state $(u_{h+1}, v_{h+1})$ through an intermediate state $(u_{h+\frac{1}{2}}, v_{h+\frac{1}{2}})$. Let $I_{1,h} = [x_{1,h}, x_{1,h+1}]$ and $I_{2,h} = [x_{2,h}, x_{2,h+1}]$ denote the spatial intervals occupied by the first and second rarefaction waves, respectively, as illustrated in Figure \ref{fig:rr}.

We apply the equivalent shock approximation to both rarefaction waves. For the first-family rarefaction connecting $(u_h, v_h)$ to $(u_{h+\frac{1}{2}}, v_{h+\frac{1}{2}})$, define:
\begin{align*}
    \xi_{u,h} &:= \frac{\int_{u_h}^{u_{h+\frac{1}{2}}} x_1(u) du}{u_{h+\frac{1}{2}} - u_h}, &
    \xi_{v,h} &:= \frac{\int_{v_h}^{v_{h+\frac{1}{2}}} x_2(v) dv}{v_{h+\frac{1}{2}} - v_h},
\end{align*}
with equivalent shock speeds $s_{u,h} = \xi_{u,h}/T$ and $s_{v,h} = \xi_{v,h}/T$.

For the second-family rarefaction connecting $(u_{h+\frac{1}{2}}, v_{h+\frac{1}{2}})$ to $(u_{h+1}, v_{h+1})$, define:
\begin{align*}
    \tilde{\xi}_{u,h} &:= \frac{\int_{u_{h+\frac{1}{2}}}^{u_{h+1}} \tilde{x}_1(u) du}{u_{h+1} - u_{h+\frac{1}{2}}}, &
    \tilde{\xi}_{v,h} &:= \frac{\int_{v_{h+\frac{1}{2}}}^{v_{h+1}} \tilde{x}_2(v) dv}{v_{h+1} - v_{h+\frac{1}{2}}},
\end{align*}
with equivalent shock speeds $\tilde{s}_{u,h} = \tilde{\xi}_{u,h}/T$ and $\tilde{s}_{v,h} = \tilde{\xi}_{v,h}/T$.

The integral relations for centered rarefaction waves in the general system \eqref{Main:eqf} yield:
\begin{equation*}
    \begin{cases}
        f_1(v_{h+\frac{1}{2}}) - f_1(v_h) = s_{u,h}(u_{h+\frac{1}{2}} - u_h), \\
        f_2(u_{h+\frac{1}{2}}) - f_2(u_h) = s_{v,h}(v_{h+\frac{1}{2}} - v_h), \\
        f_1(v_{h+1}) - f_1(v_{h+\frac{1}{2}}) = \tilde{s}_{u,h}(u_{h+1} - u_{h+\frac{1}{2}}), \\
        f_2(u_{h+1}) - f_2(u_{h+\frac{1}{2}}) = \tilde{s}_{v,h}(v_{h+1} - v_{h+\frac{1}{2}}).
    \end{cases}
\end{equation*}

Using the observable intermediate state $(u_{h+\frac{1}{2}}, v_{h+\frac{1}{2}})$ and the measured equivalent shock speeds, we obtain the reconstruction for the right state:
\begin{equation*}
    \begin{cases}
        f_1(v_{h+1}) = f_1(v_h) + s_{u,h}(u_{h+\frac{1}{2}} - u_h) + \tilde{s}_{u,h}(u_{h+1} - u_{h+\frac{1}{2}}), \\
        f_2(u_{h+1}) = f_2(u_h) + s_{v,h}(v_{h+\frac{1}{2}} - v_h) + \tilde{s}_{v,h}(v_{h+1} - v_{h+\frac{1}{2}}).
    \end{cases}
\end{equation*}

Therefore, if $f_{1,m}(v_0), \ldots, f_{1,m}(v_h)$ and $f_{2,m}(u_0), \ldots, f_{2,m}(u_h)$ are known with $f_{1,m}(v_\beta) = f_1(v_\beta)$ and $f_{2,m}(u_\alpha) = f_2(u_\alpha)$, we extend the approximations by:
\begin{equation}
    \begin{cases}
        f_{1,m}(v_{h+1}) = f_{1,m}(v_h) + s_{u,h}(u_{h+\frac{1}{2}} - u_h) + \tilde{s}_{u,h}(u_{h+1} - u_{h+\frac{1}{2}}) = f_1(v_{h+1}), \\
        f_{2,m}(u_{h+1}) = f_{2,m}(u_h) + s_{v,h}(v_{h+\frac{1}{2}} - v_h) + \tilde{s}_{v,h}(v_{h+1} - v_{h+\frac{1}{2}}) = f_2(u_{h+1}).
    \end{cases}
\end{equation}

This completes the reconstruction scheme for the case of two rarefaction waves, demonstrating that even in the absence of shock waves, the equivalent shock approximation combined with observable intermediate states provides sufficient information to reconstruct both flux functions throughout the solution domain.

\medskip

  $\textbf{(General case).}$ We now present a unified framework that extends our flux reconstruction methodology to handle arbitrary combinations of shock and rarefaction waves. This general treatment demonstrates the comprehensive applicability of our approach across all possible wave configurations that may arise in Riemann solutions.

For a fixed index $h \geq 0$, we assume that the flux approximations $f_{1,m}(v_0), \ldots, f_{1,m}(v_h)$ and $f_{2,m}(u_0), \ldots, f_{2,m}(u_h)$ are already known and satisfy $f_{1,m}(v_\alpha) = f_1(v_\alpha)$ and $f_{2,m}(u_\alpha) = f_2(u_\alpha)$. Let $\tilde{\Bu}(\cdot) = \Bu_{\mathrm{obs}}(T, \cdot)$ represent the observed solution to the Riemann problem \eqref{Main:eqf}-\eqref{eq:initial} at the fixed observation time $T$.

Due to the piecewise constant nature of our initial data and the assumption that both $f_1$ and $f_2$ possess only finite numbers of inflection points, the observed solution is guaranteed to contain only a finite number of distinct waves. However, unlike the four fundamental cases analyzed previously, the general solution may exhibit more complex wave patterns with multiple shocks and rarefactions interacting in various configurations.

To handle this general configuration systematically, we begin by identifying all wave components within the observed solution. Let $x_1 < \ldots < x_{M_1}$ denote the spatial locations of shock discontinuities in the $u$-component, and let $I_1, \ldots, I_{M_2}$ represent the intervals occupied by rarefaction waves where the solution varies continuously.

The core of our approach lies in converting each continuous rarefaction wave into an equivalent shock positioned to preserve the essential integral relations derived from the conservation laws. For a rarefaction in the $u$-component connecting states $u_l$ to $u_r$ over an interval $I_j$, we define the equivalent shock position as the mass centroid:
\begin{equation*}
\xi_{u,j}=\frac{\int_{u_{l}}^{u_r} x_1(u) \mathrm{d} u}{u_{r}-u_{l}} \in I_j .
\end{equation*}
Similarly, for the corresponding rarefaction in the $v$-component connecting states $v_l$ to $v_r$, we define:
\begin{equation*}
\xi_{v,j}=\frac{\int_{v_{l}}^{v_r} x_2(v) \mathrm{d} v}{v_{r}-v_{l}} \in I_j .
\end{equation*}

This conversion allows us to construct piecewise constant approximations $\bar{u}(\cdot)$ and $\bar{v}(\cdot)$ whose discontinuities are located at the combined set of original shock positions and equivalent shock positions. Specifically, we define the discontinuity sets as:
\begin{align*}
\{y_1, \ldots, y_M\} &= \{x_1, \ldots, x_{M_1}, \xi_{u,1}, \ldots, \xi_{u,M_2}\}, \\
\{z_1, \ldots, z_M\} &= \{x_1, \ldots, x_{M_1}, \xi_{v,1}, \ldots, \xi_{v,M_2}\},
\end{align*}
where $M = M_1 + M_2$ represents the total number of discontinuities.

The piecewise constant functions $\bar{u}(\cdot)$ and $\bar{v}(\cdot)$ take constant values between these discontinuities. Let $w_1 < \ldots < w_{M+1}$ denote the values attained by $\bar{u}$:
\begin{equation*}
\bar{u}(x)=\left\{\begin{array}{lll}
w_1=u_h, & \text { if } & x<y_1, \\
w_\alpha, & \text { if } & y_{\alpha-1}<x<y_\alpha, \alpha=2, \ldots, M, \\
w_{M+1}=u_{h+1}, & \text { if } & x>y_M .
\end{array}\right.
\end{equation*}
Similarly, let $\mu_1 < \ldots < \mu_{M+1}$ represent the values attained by $\bar{v}$:
\begin{equation*}
\bar{v}(x)=\left\{\begin{array}{lll}
\mu_1=v_h, & \text { if } & x<z_1, \\
\mu_\alpha, & \text { if } & z_{\alpha-1}<x<z_\alpha, \alpha=2, \ldots, M, \\
\mu_{M+1}=v_{h+1}, & \text { if } & x>z_M .
\end{array}\right.
\end{equation*}

The flux reconstruction in this general case follows from systematically applying the conservation laws across each discontinuity. For each equivalent shock replacing a rarefaction wave, the flux relations are given by:
\begin{equation*}
\left\{\begin{array}{l}
f_{1}(v_{l})-f_{1}(v_{r}) = s_{u,h}(u_{r}-u_{l}), \\
f_{2}(u_{r}) - f_{2}(u_{l}) = s_{v,h}(v_{r}-v_{l}),
\end{array}\right.
\end{equation*}
where $s_{u,h}=\frac{\xi_{u,j}}{T}$ and $s_{v,h}=\frac{\xi_{v,j}}{T}$ represent the equivalent shock speeds.

The complete flux reconstruction accumulates the flux changes across all discontinuities through the formulas:
\begin{equation*}
\left\{\begin{array}{l}
f_{1,m}(v_{h+1}) = f_{1,m}(v_{h}) + \sum_{\alpha=1}^{M} \frac{y_\alpha}{T}(\mu_{\alpha+1}-\mu_\alpha),\\
f_{2,m}(u_{h+1}) = f_{2,m}(u_{h}) + \sum_{\alpha=1}^{M} \frac{z_\alpha}{T}(w_{\alpha+1}-w_\alpha).
\end{array}\right.
\end{equation*}

We rigorously establish the consistency of this reconstruction by showing that:
\begin{equation*}
\left\{\begin{array}{l}
f_{1,m}(v_{h+1}) = f_{1}(v_{h+1}),\\
f_{2,m}(u_{h+1}) = f_{2}(u_{h+1}).
\end{array}\right.
\end{equation*}

This consistency follows from the fundamental property that our equivalent shock approximation preserves the integral relations derived from the conservation laws. The piecewise constant approximations $\bar{u}(\cdot)$ and $\bar{v}(\cdot)$ are constructed to satisfy the same net flux changes as the original solution $\tilde{\Bu}(\cdot)$ when integrated against the conservation laws.

\begin{remark}
    \textbf{(Observability Remark).} 
The computation of the equivalent shock positions relies only on observable quantities: the integrals $\int_{x_h}^{x_{h+1}} \tilde{u}(x) dx$ and the endpoint values $\tilde{u}(x_h), \tilde{u}(x_{h+1})$. Since $\tilde{u}$ is monotonic on each rarefaction interval, these determine the required integrals in $u$-space without explicitly computing the inverse function.
\end{remark}

\subsection{Flux Reconstruction and Convergence Analysis}

With the interpolation nodes $\{u_h\}_{h=0}^{2^m}$, $\{v_h\}_{h=0}^{2^m}$ and flux values $f_1(v_h)$, $f_2(u_h)$ determined through the analysis of Riemann solutions, we now construct the piecewise quadratic $C^1$ interpolants $f_{1,m}(v)$ and $f_{2,m}(u)$, and analyze their approximation properties. 

Given the nodes $u_* = u_0 < u_1 < \cdots < u_{2^m} = u^*$ with values $f_2(u_0), f_2(u_1), \ldots, f_2(u_{2^m})$, and $v_* = v_0 < v_1 < \cdots < v_{2^m} = v^*$ with values $f_1(v_0), f_1(v_1), \ldots, f_1(v_{2^m})$, we define the interpolants as follows.

For $f_{2,m}(u)$ on $[u_*, u^*]$:
\begin{equation*}
f_{2,m}(u) = \frac{d_{h+1}^u - d_h^u}{2\delta}(u - u_h)^2 + d_h^u(u - u_h) + f_2(u_h), \quad u \in [u_h, u_{h+1}],
\end{equation*}
where the derivatives $d_h^u = f_{2,m}'(u_h)$ are determined by the recurrence:
\begin{equation*}
d_{h+1}^u = \frac{2}{\delta}[f_2(u_{h+1}) - f_2(u_h)] - d_h^u \quad \text{for } h = 0, 1, \ldots, 2^m - 1
\end{equation*}
with initial condition:
\begin{equation*}
d_0^u = \frac{f_2(u_1) - f_2(u_0)}{\delta}.
\end{equation*}

The construction for $f_{1,m}(v)$ on $[v_*, v^*]$ follows an analogous pattern, with the mesh size $\eta$ replacing $\delta$ and the roles of $f_1$ and $f_2$ interchanged. This construction ensures $f_{2,m}(u_h) = f_2(u_h)$ and $f_{1,m}(v_h) = f_1(v_h)$ for all $h$, with $f_{2,m} \in C^1([u_*, u^*])$ and $f_{1,m} \in C^1([v_*, v^*])$.

We next analyze precise error estimates when $f_1 \in C^{1,1}[v_*, v^*]$ and $f_2 \in C^{1,1}[u_*, u^*]$ with Lipschitz constants $L_1$ for $f_1'$ and $L_2$ for $f_2'$. The analysis for both flux functions follows similar patterns, so we first present the detailed derivation for $f_{2,m}$ and then state the analogous results for $f_{1,m}$.

\textbf{Step 1: Bounding the Derivative Deviations}

For the $u$-grid with spacing $\delta = u_{h+1} - u_h$, define the average slope $A_h^u = \frac{f_2(u_{h+1}) - f_2(u_h)}{\delta}$ and deviation $\Delta_h^u = d_h^u - A_h^u$, where $d_h^u = f_{2,m}'(u_h)$ follows the recurrence:
\[
d_{h+1}^u = \frac{2}{\delta}[f_2(u_{h+1}) - f_2(u_h)] - d_h^u \quad \text{with } d_0^u = A_0^u.
\]
This recurrence translates to:
\[
\Delta_{h+1}^u = -\Delta_h^u + (A_h^u - A_{h+1}^u) \quad \text{with } \Delta_0^u = 0.
\]

To estimate $A_h^u - A_{h+1}^u$, we use the integral representation:
\[
A_h^u - A_{h+1}^u = \frac{1}{\delta} \int_{u_h}^{u_{h+1}} [f_2'(s) - f_2'(s + \delta)] ds.
\]
By Lipschitz continuity of $f_2'$ with constant $L_2$, we have $|f_2'(s) - f_2'(s + \delta)| \leq L_2\delta$ for all $s$, hence:
\[
|A_h^u - A_{h+1}^u| \leq \frac{1}{\delta} \int_{u_h}^{u_{h+1}} L_2\delta  ds = L_2\delta.
\]

The explicit solution of the recurrence is:
\[
\Delta_h^u = \sum_{j=0}^{h-1} (-1)^{h-1-j}(A_j^u - A_{j+1}^u).
\]
The alternating signs induce cancellation effects that bound the growth of the partial sums. In the worst case where the differences alternate in sign with maximum magnitude, we obtain:
\[
|\Delta_h^u| \leq 2\max_{0 \leq j \leq h-1} |A_j^u - A_{j+1}^u| \leq 2L_2\delta.
\]

By symmetry, for the $v$-grid with spacing $\eta = v_{h+1} - v_h$, we define $A_h^v = \frac{f_1(v_{h+1}) - f_1(v_h)}{\eta}$ and $\Delta_h^v = d_h^v - A_h^v$, and obtain the analogous bound:
\[
|\Delta_h^v| \leq 2L_1\eta.
\]

Thus we establish the uniform bounds:
\begin{equation*}
|\Delta_h^u| \leq 2L_2\delta, \quad |\Delta_h^v| \leq 2L_1\eta \quad \text{for all } h.
\end{equation*}

\textbf{Step 2: Function Value Errors}

We now analyze the error between the true flux functions and their piecewise quadratic approximations. For $f_2$ and its approximation $f_{2,m}$, let $\pi_2(u) = f_2(u_h) + A_h^u(u - u_h)$ be the linear interpolant of $f_2$ on the interval $[u_h, u_{h+1}]$. For any $u \in [u_h, u_{h+1}]$, we decompose the error as:

\[
|f_2(u) - f_{2,m}(u)| \leq |f_2(u) - \pi_2(u)| + |\pi_2(u) - f_{2,m}(u)|.
\]

The first term represents the linear interpolation error. By the fundamental theorem of calculus:
\[
|f_2(u) - \pi_2(u)| = \left| \int_{u_h}^u [f_2'(s) - A_h^u] ds \right|.
\]
Since $A_h^u = \frac{f_2(u_{h+1}) - f_2(u_h)}{\delta} = f_2'(\xi_h^u)$ for some $\xi_h^u \in (u_h, u_{h+1})$ by the mean value theorem, and $f_2'$ is Lipschitz continuous with constant $L_2$, we have:
\[
|f_2'(s) - A_h^u| = |f_2'(s) - f_2'(\xi_h^u)| \leq L_2|s - \xi_h^u| \leq L_2\delta.
\]
Integrating this bound over $[u_h, u]$ yields:
\[
|f_2(u) - \pi_2(u)| \leq \int_{u_h}^u L_2\delta ds \leq L_2\delta^2.
\]

The second term $R_2(u) = \pi_2(u) - f_{2,m}(u)$ arises from the quadratic correction. From the interpolation formula:
\[
f_{2,m}(u) = \frac{d_{h+1}^u - d_h^u}{2\delta}(u - u_h)^2 + d_h^u(u - u_h) + f_2(u_h),
\]
and using $d_{h+1}^u - d_h^u = -2\Delta_h^u$ from the recurrence relation, we obtain:
\[
R_2(u) = -\frac{\Delta_h^u}{\delta}(u - u_h)^2 + \Delta_h^u(u - u_h) = -\frac{\Delta_h^u}{\delta}(u - u_h)(u - u_{h+1}).
\]
The quadratic $(u - u_h)(u - u_{h+1})$ is bounded by $\frac{\delta^2}{4}$ on $[u_h, u_{h+1}]$, and $|\Delta_h^u| \leq 2L_2\delta$, giving:
\[
|R_2(u)| \leq \frac{2L_2\delta}{\delta} \cdot \frac{\delta^2}{4} = \frac{L_2}{2}\delta^2.
\]

Combining both error components:
\[
|f_2(u) - f_{2,m}(u)| \leq \frac{L_2}{2}\delta^2 + \frac{L_2}{2}\delta^2 = L_2\delta^2.
\]

By symmetry, for $f_1$ and its approximation $f_{1,m}$ on the $v$-grid with spacing $\eta$, the same analysis applies. The linear interpolation error is bounded by $\frac{L_1}{2}\eta^2$, and the quadratic residual by $\frac{L_1}{2}\eta^2$, yielding:
\[
|f_1(v) - f_{1,m}(v)| \leq L_1\eta^2.
\]

Taking suprema over the respective intervals, we establish the uniform function error estimates:
\begin{equation}
\begin{split}
\|f_2 - f_{2,m}\|_{L^\infty([u_*,u^*])} &\leq L_2\delta^2, \\
\|f_1 - f_{1,m}\|_{L^\infty([v_*,v^*])} &\leq L_1\eta^2.
\end{split}
\end{equation}

\textbf{Step 3: Derivative Errors}

We now analyze the derivative errors for the flux approximations. For $f_{2,m}$ on $[u_h, u_{h+1}]$, the derivative is:
\[
f_{2,m}'(u) = \frac{d_{h+1}^u - d_h^u}{\delta}(u - u_h) + d_h^u.
\]
Using $d_{h+1}^u - d_h^u = -2\Delta_h^u$ from the recurrence relation, this becomes:
\[
f_{2,m}'(u) = -\frac{2\Delta_h^u}{\delta}(u - u_h) + (A_h^u + \Delta_h^u) = A_h^u + \Delta_h^u \left(1 - \frac{2(u - u_h)}{\delta} \right).
\]

The derivative error decomposes as:
\[
|f_2'(u) - f_{2,m}'(u)| \leq |f_2'(u) - A_h^u| + |A_h^u - f_{2,m}'(u)|.
\]

For the first term, since $A_h^u = f_2'(\xi_h^u)$ for some $\xi_h^u \in (u_h, u_{h+1})$ and $f_2'$ is Lipschitz continuous:
\[
|f_2'(u) - A_h^u| \leq L_2|u - \xi_h^u| \leq L_2\delta.
\]

For the second term, using the expression for $f_{2,m}'(u)$:
\[
|A_h^u - f_{2,m}'(u)| = \left| \Delta_h^u \left(1 - \frac{2(u - u_h)}{\delta} \right) \right| \leq |\Delta_h^u| \cdot \max_{\theta \in [0,1]} |1 - 2\theta| = |\Delta_h^u|.
\]
With $|\Delta_h^u| \leq 2L_2\delta$, we obtain:
\[
|A_h^u - f_{2,m}'(u)| \leq 2L_2\delta.
\]

Combining both components:
\[
|f_2'(u) - f_{2,m}'(u)| \leq L_2\delta + 2L_2\delta = 3L_2\delta.
\]

By symmetry, for $f_1$ and its approximation $f_{1,m}$, the same analysis yields:
\[
|f_1'(v) - f_{1,m}'(v)| \leq 3L_1\eta.
\]

Taking suprema over the respective intervals, we establish the derivative error estimates:
\begin{equation}\label{eq:error-d}
\begin{split}
\|f_2' - f_{2,m}'\|_{L^\infty([u_*,u^*])} &\leq 3L_2\delta, \\
\|f_1' - f_{1,m}'\|_{L^\infty([v_*,v^*])} &\leq 3L_1\eta.
\end{split}
\end{equation}

\textbf{Step 4: Stability Estimate for the Coupled System}

Using the stability theory for hyperbolic conservation laws \cite{BSCR}, and noting that the characteristic speeds depend on both $f_1'$ and $f_2'$, we obtain the estimate for solution operators:
\begin{equation*}
\|S_T^{(f_1,f_2)} \Bu - S_T^{(f_{1,m},f_{2,m})} \Bu \|_{L^1} \leq C\left(\|f_1' - f_{1,m}'\|_{L^\infty} + \|f_2' - f_{2,m}'\|_{L^\infty}\right)T,
\end{equation*}
where $C$ is a positive constant independent of $\delta$ and $\eta$.

Substituting the derivative error bounds \eqref{eq:error-d} gives:
\begin{equation*}
\|S_T^{(f_1,f_2)} \Bu - S_T^{(f_{1,m},f_{2,m})} \Bu \|_{L^1} \leq 3C(L_1\eta + L_2\delta)T.
\end{equation*}

This establishes that the reconstructed flux functions $(f_{1,m}, f_{2,m})$ produce solutions that converge to the true solutions as the discretization parameters $\delta$ and $\eta$ approach zero.

\subsection{Enhanced Convergence under Higher Regularity}

We now establish improved convergence rates when the flux functions possess $C^3$ regularity. The proof strategy maintains the same framework as the $C^{1,1}$ case, but leverages the additional smoothness to obtain sharper estimates throughout the analysis.

Under the assumption that $f_1 \in C^3[v_*, v^*]$ and $f_2 \in C^3[u_*, u^*]$ with bounded third derivatives, we revisit each step of the error analysis:

\textbf{Improved bounds for deviation terms:}
Under $C^3$ regularity, the estimation of $A_h^u - A_{h+1}^u$ benefits from the additional smoothness. By Taylor expansion:
\[
A_h^u - A_h^{u+1} = \frac{1}{\delta} \int_{u_h}^{u_{h+1}} [f_2'(s) - f_2'(s + \delta)] ds,
\]
and using the bound $|f_2'(s) - f_2'(s + \delta)| \leq \|f_2'''\|_{L^\infty}\delta^2$, we obtain:
\[
|A_h^u - A_{h+1}^u| \leq \|f_2'''\|_{L^\infty}\delta^2.
\]

The alternating structure of the recurrence then gives:
\[
|\Delta_h^u| \leq 2\|f_2'''\|_{L^\infty}\delta^2,
\]
providing a quadratic improvement over the $2L_2\delta$ bound in the $C^{1,1}$ case.

\textbf{Enhanced function value errors:}
For the linear interpolation error, expand \(f_2(u)\) and its linear interpolant \(\pi_2(u)\) around \(u_h\) for \(u = u_h + t\delta\) (\(t \in [0,1]\)):
\[
f_2(u) = f_2(u_h) + t\delta f_2'(u_h) + \frac{(t\delta)^2}{2}f_2''(u_h) + \frac{(t\delta)^3}{6}f_2'''(a), \quad a \in (u_h, u),
\]
\[
\pi_2(u) = f_2(u_h) + t\left[\delta f_2'(u_h) + \frac{\delta^2}{2}f_2''(u_h) + \frac{\delta^3}{6}f_2'''(b)\right], \quad b \in (u_h, u_{h+1}).
\]
Subtracting these gives the error:
\[
f_2(u) - \pi_2(u) = \frac{\delta^2}{2}f_2''(u_h)(t^2 - t) + \frac{\delta^3}{6}\left[t^3f_2'''(a) - tf_2'''(b)\right].
\]
Bounding using \(|t^2 - t| \leq 1/4\) (since \(t \in [0,1]\)), \(|t^3|, |t| \leq 1\), and \(|f_2''(u_h)| \leq \|f_2'''\|_{L^\infty}\delta\) (by \(C^3\) regularity, as \(f_2''\) is Lipschitz continuous with constant \(\|f_2'''\|_{L^\infty}\)):

- The first term contributes: \(\frac{\delta^2}{2} \cdot \|f_2'''\|_{L^\infty}\delta \cdot \frac{1}{4} = \frac{1}{8}\|f_2'''\|_{L^\infty}\delta^3\)

- The second term contributes: \(\frac{\delta^3}{6} \cdot \left(1 + 1\right)\|f_2'''\|_{L^\infty} = \frac{1}{3}\|f_2'''\|_{L^\infty}\delta^3\)

Summing these gives:
\[
|f_2(u) - \pi_2(u)| \leq \left(\frac{1}{8} + \frac{1}{3}\right)\|f_2'''\|_{L^\infty}\delta^3 = \frac{11}{24}\|f_2'''\|_{L^\infty}\delta^3.
\]

For the residual \(R_2(u) = \pi_2(u) - f_{2,m}(u)\), recall \(R_2(u) = -\frac{\Delta_h^u}{\delta}(u - u_h)(u - u_{h+1})\). Using \(|\Delta_h^u| \leq 2\|f_2'''\|_{L^\infty}\delta^2\) and \(\max|(u - u_h)(u - u_{h+1})| = \delta^2/4\) (attained at the midpoint of \([u_h, u_{h+1}]\)):
\[
|R_2(u)| \leq \frac{2\|f_2'''\|_{L^\infty}\delta^2}{\delta} \cdot \frac{\delta^2}{4} = \frac{1}{2}\|f_2'''\|_{L^\infty}\delta^3.
\]

Combining these gives cubic convergence:
\[
\|f_2 - f_{2,m}\|_{L^\infty} \leq \left(\frac{11}{24} + \frac{1}{2}\right)\|f_2'''\|_{L^\infty}\delta^3 = \frac{23}{24}\|f_2'''\|_{L^\infty}\delta^3 \leq \|f_2'''\|_{L^\infty}\delta^3.
\]

\textbf{Improved derivative approximations:}
For derivative errors, both components exhibit enhanced convergence under $C^3$ regularity:
\begin{itemize}
    \item The error between $f_2'$ and the linear interpolant's derivative satisfies:
    \[
    |f_2'(u) - A_h^u| \leq \frac{1}{2}\|f_2'''\|_{L^\infty}\delta^2
    \]
    \item The deviation term satisfies:
    \[
    |A_h^u - f_{2,m}'(u)| \leq |\Delta_h^u| \leq 2\|f_2'''\|_{L^\infty}\delta^2
    \]
\end{itemize}
Combining these gives the quadratic convergence:
\[
\|f_2' - f_{2,m}'\|_{L^\infty} \leq \frac{5}{2}\|f_2'''\|_{L^\infty}\delta^2.
\]

By symmetry, identical improvements apply to $f_1$ and $f_{1,m}$, giving:
\[
\|f_1 - f_{1,m}\|_{L^\infty} \leq \|f_1'''\|_{L^\infty}\eta^3, \quad \|f_1' - f_{1,m}'\|_{L^\infty} \leq \frac{5}{2}\|f_1'''\|_{L^\infty}\eta^2.
\]

The corresponding improvement in solution error follows from the stability estimate:
\[
\left\|S_T^{(f_1,f_2)} \Bu - S_T^{(f_{1,m},f_{2,m})} \Bu\right\|_{\mathbf{L}^1} \leq C T (\|f_1'''\|_{L^\infty}\eta^2 + \|f_2'''\|_{L^\infty}\delta^2),
\]
where $C$ is the same stability constant as in the $C^{1,1}$ case.

This remarkable improvement demonstrates that our reconstruction method benefits substantially from higher regularity of the unknown flux functions. In practical applications where physical fluxes are typically smooth, the method achieves rapid convergence, making it particularly effective for high-precision flux identification in continuum mechanics and related fields.

\section{Application to Isentropic Euler Equations and the p-System}\label{Sec:Application}

The isentropic Euler equations represent a cornerstone model in compressible fluid dynamics, describing the motion of inviscid fluids under adiabatic conditions, see \cite{DiPerna1983, Lions1998, Chen1999}. While the full Euler system comprises three conservation laws (mass, momentum, and energy), the isentropic approximation reduces this to a 2×2 system under the assumption of constant entropy. This is achieved by replacing the energy conservation law with an isentropic equation of state $p = p(\rho)$, which provides a thermodynamic closure while eliminating the need for explicit energy tracking. The isentropic assumption is physically justified in scenarios where heat transfer and dissipative effects are negligible, such as in rapid compression/expansion processes where the flow remains effectively adiabatic and reversible. This transformation makes the isentropic Euler equations an ideal testbed for our flux inversion methodology, demonstrating its applicability to physically fundamental systems while maintaining the mathematical structure required by our theoretical framework.

The isentropic Euler equations form a 2×2 system expressed as:
\begin{equation}
    \begin{cases}
        \dfrac{\partial \rho}{\partial t} + \dfrac{\partial (\rho u)}{\partial x} = 0, \\[10pt]
        \dfrac{\partial (\rho u)}{\partial t} + \dfrac{\partial}{\partial x} \left( \rho u^2 + p(\rho) \right) = 0,
    \end{cases}
\end{equation}
where $\rho$ denotes the fluid density, $u$ the velocity, and $p(\rho)$ the pressure determined by the equation of state. This system can be cast into our general framework by defining the conserved variables and fluxes as:
\begin{equation}
    \Bu = \begin{pmatrix} \rho \\ \rho u \end{pmatrix}, \quad F(\Bu) = \begin{pmatrix} \rho u \\ \rho u^2 + p(\rho) \end{pmatrix}.
\end{equation}

\begin{remark}[Lagrangian Formulation and the p-System]
The isentropic Euler equations admit an important Lagrangian formulation through the coordinate transformation:
\begin{equation}
    \begin{cases}
        s = t,\\
        y = \int_{(0,0)}^{(x,t)} \rho dx - \rho u dt.
    \end{cases}
\end{equation}
Under this transformation, the system reduces to the p-system:
\begin{equation}
    \begin{cases}
        \partial_s v - \partial_y u = 0,\\
        \partial_s u + \partial_y p = 0,
    \end{cases}
\end{equation}
where the specific volume $v = 1/\rho$ replaces density as the primary variable. According to \cite{Wagner1987}, this transformation establishes the equivalence of solutions between the isentropic Euler equations and the p-system. This mathematical equivalence demonstrates that our flux inversion methodology applies to both formulations, with the p-system being particularly useful for problems involving material interfaces and providing an alternative perspective for analyzing wave propagation in compressible media.
\end{remark}

The reduction from the full Euler system to this $2\times 2$ isentropic system is physically justified when heat transfer and dissipative effects are negligible, a scenario commonly encountered in high-speed aerodynamics, gas dynamics, and certain astrophysical applications. This reduction is crucial for the application of our method, as it preserves the mathematical structure for which our wave front tracking and flux reconstruction techniques were developed.

The characteristic structure of this system reveals why our flux inversion approach is particularly powerful in this context. For the Eulerian coordinate formulation, the Jacobian matrix
\begin{equation}
    DF(\Bu) = \begin{pmatrix}
        0 & 1 \\
        -u^2 + p'(\rho) & 2u
    \end{pmatrix}
\end{equation}
has eigenvalues $\lambda_{1,2} = u \mp \sqrt{p'(\rho)}$, which correspond to the speeds of acoustic waves propagating upstream and downstream relative to the fluid motion. For the Lagrangian p-system formulation, the characteristic speeds are $\lambda_{1,2} = \mp\sqrt{-p'(v)}$ with $v = 1/\rho$. In both cases, the observable wave patterns in Riemann solutions-whether shocks or rarefaction waves-directly encode information about the unknown pressure function through these characteristic speeds and the corresponding Rankine-Hugoniot conditions.

The application of Theorem \ref{main:thm-f} to both formulations of compressible flow yields the following significant result:

\begin{thm}\label{thm:euler}
For the isentropic Euler equations (and equivalently for the p-system) with pressure function $p(\rho) \in C^{1,1}[\rho_*, \rho^*]$ having Lipschitz constant $L_p$ for its derivative, there exists a piecewise quadratic $C^1$ approximation $p_m(\rho)$ such that:
\begin{equation}
\|p - p_m\|_{L^\infty([\rho_*,\rho^*])} \leq L_p\delta^2, \quad \|p' - p_m'\|_{L^\infty([\rho_*,\rho^*])} \leq 3L_p\delta,
\end{equation}
where $\delta$ is the discretization parameter. Moreover, the reconstructed pressure function produces flow solutions that approximate the true dynamics in both formulations:
\begin{equation}\label{Stability-P}
\left\|S_{T}^{p_m} \Bu - S_{T}^{p}\Bu\right\|_{\mathbf{L}^1} \leqslant C T \delta
\end{equation}
for all $\Bu \in \mathscr{D}^F$ with values in $[\rho_*, \rho^*] \times \mathbb{R}$, where $C$ is a constant independent of $\delta$.
\end{thm}

The application of our flux inversion method to these equivalent systems enables the identification of the equation of state from dynamic flow observations. This has profound implications for experimental fluid dynamics, where determining the correct equation of state traditionally requires separate thermodynamic measurements \cite{Thompson1972}. Our approach instead extracts this fundamental relation directly from wave propagation behavior. Consider a shock tube experiment designed to generate Riemann-type initial conditions with states $(\rho_L, u_L)$ and $(\rho_R, u_R)$. Upon diaphragm rupture, the resulting wave structure-typically a shock wave propagating into one region and a rarefaction wave expanding into the other-encodes the essential information about the pressure function $p(\rho)$. Through a sequence of such experiments with systematically varied initial conditions, and by observing the wave positions and speeds, our method reconstructs $p(\rho)$ over the density range $[\rho_*, \rho^*]$. Equivalently, similar wave patterns in the Lagrangian p-system framework \cite{Smoller1994} can be utilized, demonstrating the method's flexibility.

The observable wave structures-whether shocks or rarefactions-provide distinct information about the pressure function through the Rankine-Hugoniot conditions for discontinuities and the integral curves for smooth regions. Our equivalent shock concept proves particularly valuable for handling continuous rarefaction waves, allowing a unified treatment within the reconstruction framework that bridges discontinuous and continuous solution features.

The significance of this application extends to practical engineering challenges, such as in high-speed aerodynamics where an accurate equation of state is crucial for predicting shock positions and overall flow behavior. Our method provides a pathway to validate or discover equations of state for novel working fluids under extreme conditions. Furthermore, the convergence guarantees established in Theorem \ref{thm:euler} ensure that as experimental resolution improves, the reconstructed pressure function converges to the true physical relation, providing mathematical confidence for critical engineering applications.

Beyond experimental fluid dynamics, the stability estimate \eqref{Stability-P} provides crucial theoretical foundations for computational methods. While modern high-resolution schemes-including Godunov-type methods \cite{LeVeque2002, Toro1999}, approximate Riemann solvers \cite{Roe1981}, and discontinuous Galerkin methods \cite{Cockburn1998} have demonstrated empirical success, our results provide a theoretical foundation for their convergence by establishing rigorous error estimates for flux function approximations. The estimate quantifies how flux approximation errors propagate to solution errors, addressing a fundamental question in the numerical analysis of conservation laws. This is particularly relevant for the p-system and isentropic Euler equations, where solvers like Roe's method \cite{Roe1981} rely on local linearization of the flux Jacobian. Our analysis provides a mathematical framework for understanding the accumulation of such approximation errors, thereby complementing numerical experimentation and offering theoretical guidance for adaptive mesh refinement and error control.

The successful application of our flux inversion methodology thus demonstrates its capacity to address fundamental problems in fluid dynamics, bridging the gap between abstract mathematical theory and practical applications while opening new avenues for both experimental and computational advancements.

\vspace{0.4cm}
\noindent\textbf{Acknowledgment.} 
 The work of H. Liu is supported by the Hong Kong RGC General Research Funds (projects 11311122, 11300821, and 11303125), the NSFC/RGC Joint Research Fund (project  N\_CityU101/21), the France-Hong Kong ANR/RGC Joint Research Grant, A-CityU203/19. The work of Y. Jiang is supported by the China Natural National Science Foundation (No. 123B2017).

\newpage
\renewcommand\refname{References}

\end{document}